\documentclass{gtmon_a}
\pdfoutput=1

\usepackage{amscd}

\proceedingstitle{Groups, homotopy and configuration spaces (Tokyo
  2005)}
\conferencestart{5 July 2005}
\conferenceend{11 July 2005}
\conferencename{Groups, homotopy and configuration spaces,
                in honour of Fred Cohen's 60th birthday}
\conferencelocation{University of Tokyo, Japan}

\editor{Norio Iwase}
\givenname{Norio}
\surname{Iwase}

\editor{Toshitake Kohno}
\givenname{Toshitake}
\surname{Kohno}

\editor{Ran Levi}
\givenname{Ran}
\surname{Levi}

\editor{Dai Tamaki}
\givenname{Dai}
\surname{Tamaki}

\editor{Jie Wu}
\givenname{Jie}
\surname{Wu}

\title[Lie algebras of symplectic derivations and cycles on the moduli spaces]{Lie algebras of symplectic derivations\\and cycles on the moduli spaces}
\author{Shigeyuki Morita}
\givenname{Shigeyuki}
\surname{Morita}
\address{Graduate School of Mathematical Sciences\\
University of Tokyo\\\newline
Komaba, Tokyo 153-8914\\Japan}
\email{morita\char'100ms.u-tokyo.ac.jp}
\urladdr{}

\volumenumber{13}
\issuenumber{}
\publicationyear{2008}
\papernumber{16}
\startpage{335}
\endpage{354}

\doi{}
\MR{}
\Zbl{}

\arxivreference{math/0608673}


\received{21 August 2006}
\revised{28 March 2007}
\accepted{30 March 2007}
\published{25 February 2008}
\publishedonline{25 February 2008}
\proposed{}
\seconded{}
\corresponding{}
\version{}

\makeatletter
\def\cnewtheorem#1[#2]#3{\newtheorem{#1}{#3}
\expandafter\let\csname c@#1\endcsname\c@theorem}

\AtBeginDocument{\let\bar\wbar\let\tilde\wtilde\let\hat\what}

\newtheorem{theorem}{Theorem}%[section]
\cnewtheorem{lemma}[theorem]{Lemma}
\cnewtheorem{proposition}[theorem]{Proposition}
\cnewtheorem{corollary}[theorem]{Corollary}
\theoremstyle{definition}
\cnewtheorem{definition}[theorem]{Definition}
\cnewtheorem{problem}[theorem]{Problem}
\cnewtheorem{example}[theorem]{Example}
\cnewtheorem{conjecture}[theorem]{Conjecture}
\theoremstyle{remark}
\cnewtheorem{remark}[theorem]{Remark}

\makeatother  %  move after \newtheorem block
\def\ra{{\rightarrow}}
\def\lra{{\longrightarrow}}
\def\om{{\omega}}

\def\al{{\alpha}}

\def\ga{{\gamma}}

\def\bZ{{\mathbb Z}}
\def\bQ{{\mathbb Q}}

\def\M1g{{\mathcal M}_{g,1}}
\def\I1g{{\mathcal I}_{g,1}}

\newcommand\Hom{\operatorname{Hom}}

\newcommand\sgn{\operatorname{sgn}}

\newcommand\Ker{\operatorname{Ker}}

\begin{document}

\begin{abstract}   % type your abstract below
We consider the Lie algebra consisting of all derivations on the
free associative algebra, generated by the first homology group
of a closed oriented surface, which kill the symplectic class.
We find the first non-trivial abelianization of
this Lie algebra and discuss its relation to unstable cohomology
classes of the moduli space of curves via a theorem of Kontsevich.

\end{abstract}
\maketitle

\section{Introduction}

In \cite{Kontsevich93,Kontsevich94}, Kontsevich considered three
infinite dimensional Lie algebras, denoted by ${\mathfrak c}_g,
{\mathfrak a}_g, {\mathfrak l}_g$, and described the stable homology
of them, as $g$ tends to $\infty$, in terms of his {\it graph
homology} for the case ${\mathfrak c}_\infty$, the rational cohomology
of the moduli spaces of curves with unlabeled punctures for the case
${\mathfrak a}_\infty$, and the rational cohomology of the outer
automorphism group ${\mathrm Out}\, F_n$ of free groups $F_n\ (n\geq
2)$ for the case ${\mathfrak l}_\infty$, respectively. See also Conant
and Vogtmann \cite{CoV03} for more details on this result.

In this paper, we consider the latter two Lie algebras. The ideal
${\mathfrak l}^+_g$ of the last Lie algebra ${\mathfrak l}_g$,
consisting of derivations with {\it positive} degrees, appeared,
before the work of Kontsevich cited above, in the study of Johnson
homomorphisms which were introduced (by Johnson) in
\cite{Johnson80,Johnson83}; see our papers \cite{Morita89a,Morita91}.
In \cite{Morita93a}, we constructed a large abelian quotient of the
Lie algebra ${\mathfrak l}^+_g$ and by making use of this, we defined
in \cite{Morita99} many homology classes of $\mathrm{Out}\, F_n$ via
the above theorem of Kontsevich.  These classes have been investigated
by Conant and Vogtmann \cite{CoV04,CoV06} by interpreting them
geometrically on the Outer Space of Culler and Vogtmann \cite{CV}.

In this paper, we consider the associative case $\mathfrak{a}_{g}$.
We show the first non-trivial abelian quotient of the ideal
$\mathfrak{a}^+_{g}$. More precisely, we determine the weight $2$
part of $H_1(\mathfrak{a}^+_{g})$ where the grading by weights is induced by
that of the graded Lie algebra $\mathfrak{a}^+_{g}$.
We discuss possible application of
this abelianization to the problem of finding unstable cohomology classes
of the moduli space of curves. We also consider a generalization of our result
to the case of the Lie algebra of derivations of a free associative
alegbra {\it without} the symplectic constraint.

\section{Two Lie algebras of symplectic derivations}

In this section, we recall the definition of the two Lie algebras
$\mathfrak{a}_{g}, \mathfrak{l}_{g}$ mentioned in the
previous section briefly.

Let $\Sigma_g$ be a closed oriented surface of genus $g\geq 1$.
We simply denote by $H$ the first integral homology group
$H_1(\Sigma_g;\bZ)$ and let $H_\bQ=H\otimes\bQ=H_1(\Sigma_g;\bQ)$.
The intersection pairing
$$
H\times H\ni (u,v)\longmapsto u\cdot v\in\bZ,
$$
which is a non-degenerate skew symmetric pairing,
induces the Poincar\'e duality isomorphisms
$$
H\cong H^*,\quad H_\bQ\cong H_\bQ^*.
$$
Let $\mathrm{Sp}(H_\bQ)$ denote the automorphism group of $H_\bQ$
preserving the above pairing. If we choose a symplectic basis of
$H_\bQ$, then $\mathrm{Sp}(H_\bQ)$ can be identified with the
algebraic group $\mathrm{Sp}(2g,\bQ)$. In short, we understand
$H_\bQ$ as a symplectic vector space of dimension $2g$.

Let ${\mathcal L}_{g}$ denote the free graded Lie algebra
generated by $H$ and let ${\mathcal L}^\bQ_{g}={\mathcal L}_{g}\otimes\bQ$.
We denote by ${\mathcal L}_{g}^\bQ(k)$ the degree $k$ part of
${\mathcal L}_{g}^\bQ$. In particular, we have
$$
{\mathcal L}_{g}^\bQ(1)=H_\bQ,\quad
{\mathcal L}_{g}^\bQ(2)\cong \Lambda^2 H_\bQ,\quad
{\mathcal L}_{g}^\bQ(3)\cong H_\bQ\otimes \Lambda^2 H_\bQ/\Lambda^3 H_\bQ
$$
and so on. Next let $T(H_\bQ)$ denote the tensor algebra
(without unit) on $H_\bQ$. In other
words, it is the free associative algebra without unit generated by $H_\bQ$.
It is a graded algebra whose degree $k$ part is
$H_\bQ^{\otimes k}$.
As is well known, we can consider ${\mathcal L}_{g}^\bQ$ as
a natural submodule of $T(H_\bQ)$ consisting of Lie elements.
We denote by
$$
\omega_0\in {\mathcal L}_{g}^\bQ(2)\cong\Lambda^2 H_\bQ
\subset H_\bQ^{\otimes 2}
$$
the symplectic class.

\begin{definition}
We define two Lie algebras ${\mathfrak a}_g$ and ${\mathfrak l}_g$ by setting
\begin{align*}
{\mathfrak a}_g&=\{\mathrm{derivation}\ D\ \mathrm{of}\ \mathrm{the}\
\mathrm{tensor}\ \mathrm{algebra}\
T(H_\bQ); D(\om_0)=0\}\\
{\mathfrak l}_g&=\{\mathrm{derivation}\ D\ \mathrm{of}\
\mathrm{the}\ \mathrm{free}\ \mathrm{Lie}\ \mathrm{algebra}\
{\mathcal L}_{g}^\bQ:D(\om_0)=0\}.
\end{align*}
\end{definition}

These Lie algebras are naturally graded by the degrees of derivations.
More precisely, we have direct sum decompositions
$$
{\mathfrak a}_g=\bigoplus_{k=0}^\infty {\mathfrak a}_g(k),\quad
{\mathfrak l}_g=\bigoplus_{k=0}^\infty {\mathfrak l}_g(k)
$$
\begin{align*}
{\mathfrak a}_g(k)&=\{D\in \Hom \left(H_\bQ,H_\bQ^{\otimes 
(k+1)}\right);D(\om_0)=0\}\tag*{\hbox{where}} \\
{\mathfrak l}_g(k)&=\{D\in \Hom \left(H_\bQ,{\mathcal L}_{g}^\bQ(k+1)\right)
;D(\om_0)=0\}
\end{align*}
denote the degree $k$ summands.

It is easy to see that
$$
{\mathfrak a}_g(0)={\mathfrak l}_g(0)\cong\mathfrak{sp}(2g,\bQ)\cong S^2 H_\bQ
$$
where $\mathfrak{sp}(2g,\bQ)$ denotes the Lie algebra of $\mathrm{Sp}(2g,\bQ)$
and $S^2 H_\bQ$ denotes the second symmetric power of $H_\bQ$.
Also observe that ${\mathfrak l}_g(k)$ can be considered as a natural
submodule of ${\mathfrak a}_g(k)$ so that
$\mathfrak{l}_g$ is a Lie subalgebra of $\mathfrak{a}_g$.
$$
{\mathfrak a}^+_g=\bigoplus_{k=1}^\infty {\mathfrak a}_g(k),\quad
{\mathfrak l}^+_g=\bigoplus_{k=1}^\infty {\mathfrak l}_g(k)\leqno{\rm Let}
$$ denote the ideal of $\mathfrak{a}_g$ and
$\mathfrak{l}_g$,respectively, consisting of elements with positive
degrees. Then, as we already mentioned in the introduction, the latter
one $\mathfrak{l}^+_g$ is the same as the Lie algebra
$\mathfrak{h}_{g,1}^\bQ$ considered in the study of Johnson
homomorphisms (see \cite{Morita89a,Morita91}) before the work of
Kontsevich.  In this paper, however, we use Kontsevich's notation
$\mathfrak{l}^+_g$ instead of ours.

\begin{proposition}\label{prop:al}
The Poincar\'e duality $H_\bQ\cong H_\bQ^*$
induces canonical isomorphisms
\begin{align*}
{\mathfrak a}_g(k)&=\{D\in\Hom\left(H_\bQ,H_\bQ^{\otimes (k+1)}\right); 
D(\omega_0)=0\}\\
&\cong \left(H_\bQ^{\otimes (k+2)}\right)^{\bZ/(k+2)\bZ}
\end{align*}
for the associative case and
\begin{align*}
{\mathfrak l}_g(k)&=\{D\in\Hom\left(H_\bQ,{\mathcal L}_{g}^\bQ(k+1)\right)
; D(\omega_0)=0\}\\
&\cong \Ker\left(H_\bQ\otimes{\mathcal L}_{g}^\bQ(k+1)\
\overset{[\ ,\ ]}{\lra}\ {\mathcal L}_{g}^\bQ(k+2)\right)
\end{align*}
for the Lie case, where for the former case the cyclic group $\bZ/(k+2)\bZ$ 
of order $k+2$
acts naturally on $H_\bQ^{\otimes (k+2)}$ by cyclic permutations.
\end{proposition}

\begin{proof}
The latter statement was proved in \cite{Morita93a} and the former statement
was mentioned already by Kontsevich \cite{Kontsevich93,Kontsevich94}.
Here we give a proof of the former for completeness.
Choose a symplectic basis
$x_1,\cdots,x_g,y_1,\cdots,y_g$ of $H_\bQ$. Then it can be checked that the
isomorphism
$$
\Hom\left(H_\bQ,H_\bQ^{\otimes (k+1)}\right)=
H_\bQ^*\otimes H_\bQ^{\otimes (k+1)}
\cong
H_\bQ\otimes H_\bQ^{\otimes (k+1)}
$$
induced by the Poincar\'e duality $H_\bQ^*\cong H_\bQ$ is given by
the correspondence
$$
\Hom\left(H_\bQ,H_\bQ^{\otimes (k+1)}\right)\ni D\longmapsto
D^*\in H_\bQ\otimes H_\bQ^{\otimes (k+1)}
$$
where
$$
D^*=\sum_{i=1}^g \{x_i\otimes D(y_i)-y_i\otimes D(x_i)\}.
$$
Now assume that $D$ belongs to $\mathfrak{a}_g(k)$, namely
$D(\om_0)=0$. Since $D$ acts on $\om_0=\sum_i (x_i\otimes y_i
-y_i\otimes x_i)$ by a derivation, we have
\begin{align*}
D(\om_0)&=
\sum_{i=1}^g
\Big\{D(x_i)\otimes y_i +x_i\otimes D(y_i) -D(y_i)\otimes x_i
-y_i\otimes D(x_i)\Big\}\\
&=D^* - \sigma_{k+2}(D^*)
\end{align*}
where
$
\sigma_{k+2}
$
is a generator of the cyclic group $\bZ/(k+2)\bZ$
acting cyclically on
$H_\bQ\otimes H_\bQ^{\otimes (k+1)}=H_\bQ^{\otimes (k+2)}$ by
$$
\sigma_{k+2}(u_1\otimes u_2\otimes\cdots\otimes u_{k+2})
=u_2\otimes\cdots\otimes u_{k+2}\otimes u_1.
$$
It follows that $D(\om_0)=0$ if and only if $D^*$ is invariant
under the action of $\bZ/(k+2)\bZ$. This completes the proof.
\end{proof}

\begin{example}\label{ex:1}
If we apply the above proposition to the degree $1$ cases,
we see that
\begin{gather*}
{\mathfrak a}_g(1)\cong \left(H_\bQ^{\otimes 3}\right)^{\bZ/3}
\cong S^3 H_\bQ\oplus\Lambda^3 H_\bQ\\
{\mathfrak l}_g(1)\cong
\Ker\left(H_\bQ\otimes{\mathcal L}_{g}^\bQ(2)\
\overset{[\ ,\ ]}{\lra}\ {\mathcal L}_{g}^\bQ(3)\right)
\cong \Lambda^3 H_\bQ\tag*{\hbox{and}}
\end{gather*}
where $S^3 H_\bQ$ denotes the third symmetric power of $H_\bQ$.

In fact, as is well known, there is a canonical decomposition
$$
H_\bQ^{\otimes 3}\cong S^3 H_\bQ\oplus \Lambda^3 H_\bQ\oplus
V\oplus V
$$
where
$V$ denotes a certain $\mathrm{GL}(H_\bQ)$--irreducible representation.
It is easy to see that both $S^3 H_\bQ$ and $\Lambda^3 H_\bQ$
are $\bZ/3\bZ$--invariant while there is no non-trivial $\bZ/3\bZ$--invariant
subspace of $V\oplus V$. The latter statement follows from
a well known fact that there exists a natural isomorphism
$$
\mathcal{L}^\bQ_{g}(3)\cong H_\bQ\otimes \mathcal{L}^\bQ_{g}(2)/
\Lambda^3 H_\bQ
$$
where the factor $\Lambda^3 H_\bQ$ corresponds precisely to the
Jacobi identity of the free Lie algebra $\mathcal{L}^\bQ_{g}$.
\end{example}

\section{Main result}

In our paper \cite{Morita93a}, we proved the following theorem.

\begin{theorem}{\rm\cite{Morita93a}}\qua
There exists a surjective homomorphism
$$
{\mathfrak l}_g^+\lra \Lambda^3 H_\bQ\oplus
\bigoplus_{k=1}^\infty S^{2k+1} H_\bQ
$$
of graded Lie algebras where $S^{2k+1} H_\bQ$ denotes the
$(2k+1)$-st symmetric power of $H_\bQ$ and the target is considered
as an {\it abelian} Lie algebra. In particular, the Lie algebra
${\mathfrak l}_g^+$ is not finitely generated.
\end{theorem}

In view of this result, we have the following important problem.

\begin{problem}
Determine whether the Lie algebra ${\mathfrak a}^+_g$ is finitely generated
or not. In particular, is the abelianization
$H_1({\mathfrak a}^+_g)$ of ${\mathfrak a}^+_g$ finitely generated or not?
\end{problem}

As the first step towards the solution of the above problem,
we determined the weight $2$ part of the abelianization
of ${\mathfrak a}^+_g$ as follows.
Recall here that the Lie algebra ${\mathfrak a}^+_g$ is graded
so that its homology group $H_*({\mathfrak a}^+_g)$ is bigraded.
In particular, we have a direct sum decomposition
$$
H_1({\mathfrak a}^+_g)=\bigoplus_{m=1}^\infty
H_1({\mathfrak a}^+_g)_m
$$
where $H_1({\mathfrak a}^+_g)_m$ denotes the subspace generated by
$1$--cycles which are homogeneous of degree $m$. More precisely
$$
H_1({\mathfrak a}^+_g)_m=\mathrm{quotient}\ \mathrm{of}\
\mathfrak{a}_g(m)\ \mathrm{by}\
\sum_{i+j=m,i,j>0} [\mathfrak{a}_g(i),\mathfrak{a}_g(j)].
$$
In particular, we have
$$
H_1({\mathfrak a}^+_g)_1=\mathfrak{a}_g(1)\cong S^3H_\bQ\oplus \Lambda^3H_\bQ.
$$

The following is the main theorem of this paper.

\begin{theorem}\label{th:main}
We have an isomorphism
$$
H_1({\mathfrak a}_g^+)_{2}\cong \Lambda^2 H_\bQ/\bQ(\omega_0)
$$
of $\mathrm{Sp}(H_\bQ)$--modules. More precisely,
the correspondence
\begin{equation}\label{eq:h}
H_\bQ^{\otimes 4}\supset {\mathfrak a}_g(2)\
\overset{C_{13}}{\lra}\
H_\bQ^{\otimes 2} \lra\Lambda^2 H_\bQ \ \mathrm{mod}\ \bQ(\omega_0)
\end{equation}
induces the above isomorphism. Here
$$
C_{13}:H_\bQ^{\otimes 4}\lra H_\bQ^{\otimes 2}
$$
denotes the mapping given by
$$
C_{13}(u_1\otimes u_2\otimes u_3\otimes u_4)=
(u_1\cdot u_3)\, u_2\otimes u_4\quad (u_i\in H_\bQ)
$$
and the last mapping in \eqref{eq:h} is the natural projection.
\end{theorem}

\section{Proof of the main result}

In this section, we prove \fullref{th:main}.  We begin by
determining the $\mathrm{Sp}(H_\bQ)$--irreducible decomposition of the
degree $2$ part ${\mathfrak a}_g(2)$ of the Lie algebra
$\mathfrak{a}_g$.  Here we use the terminology of \cite{Morita01} to
describe irreducible representations of the algebraic group
$\mathrm{Sp}(2g,\bQ)$ as well as their highest weight vectors (see
also Fulton and Harris \cite{FH} for the generality of the
representation theory of $\mathrm{Sp}(2g,\bQ)$). In particular, we use
the symbol $[b_1b_2\cdots b_k]$ for expressing irreducible
representations which are in one to one correspondence with Young
diagrams. For example, the $k$-th symmetric power $S^k H_\bQ$ is
denoted by the symbol $[k]$ and $[1^2]$ denotes the irreducible
representation $\Lambda^2H_\bQ/\bQ(\om_0)$.  We denote by $\bQ$ the
one dimensional trivial representation which corresponds to the empty
Young diagram.  We also fix a symplectic basis
$x_1,\cdots,x_g,y_1,\cdots,y_g$ of $H_\bQ$ which we use to describe
highest weight vectors.  For example, the highest weight vector of
$[1^2]$ is $x_1\wedge x_2$.

\begin{lemma}\label{lem:4}
The irreducible decomposition of the $\mathrm{Sp}(H_\bQ)$--module
$H_\bQ^{\otimes 4}$ is given by
$$
H_\bQ^{\otimes 4}\cong 3 \bQ\oplus 6 [1^2]\oplus 6 [2]\oplus
2 [2^2]\oplus 3 [31] \oplus 3 [21^2]\oplus [1^4]\oplus [4]
$$
for all $g\geq 4$.
\end{lemma}

\begin{proof}
This can be shown by a direct computation using the well known
decomposition formula for tensor products of
$\mathrm{Sp}(H_\bQ)$--irreducible representations (see eg \cite{FH}).
\end{proof}

\begin{proposition}\label{prop:a2}
We have a canonical decomposition
$$
{\mathfrak a}_g(2)\cong
\bQ\oplus 2 [1^2]\oplus [2]\oplus [2^2]\oplus [21^2]\oplus
[4]
$$
for all $g\geq 4$.
\end{proposition}

\begin{proof}
By \fullref{prop:al}, we have an isomorphism
$$
\mathfrak{a}_g(2)\cong\left(H_\bQ^{\otimes 4}\right)^{\bZ/4\bZ}.
$$
Hence we have only to determine the $\bZ/4\bZ$--invariant part
in the decomposition given in \fullref{lem:4} above.
We consider each irreducible component case by case.

\smallskip
{\bf Case (i)}\qua (the trivial representation $\bQ$)
\par
\noindent
We define three $\mathrm{Sp}(H_\bQ)$--invariant elements of $H_\bQ^{\otimes 4}$
by setting
$$
\om_{12}=\om_0\otimes\om_0,\quad
\om_{13}=\om_0\otimes_{13}\om_0,\quad
\om_{14}=\om_0\otimes_{14}\om_0
$$
where $\otimes_{13}$ and $\otimes_{14}$ are defined by
\begin{align*}
(u_1\otimes u_2)\otimes_{13}(u_3\otimes u_4)&=u_1\otimes u_3\otimes 
u_2\otimes u_4\\
(u_1\otimes u_2)\otimes_{14}(u_3\otimes u_4)&=u_1\otimes u_3\otimes 
u_4\otimes u_2.
\end{align*}
It is easy to see that these three elements form a basis
of the $\mathrm{Sp}(H_\bQ)$--trivial part of
$H_\bQ^{\otimes 4}$. Let $\sigma_4\in \bZ/4\bZ$ be the generator as before.
Then it acts on the above elements as
$$
\sigma_4(\om_{12})=-\om_{14},\quad
\sigma_4(\om_{13})=-\om_{13},\quad
\sigma_4(\om_{14})=-\om_{12}.
$$
It follows that the dimension of the $\mathrm{Sp}(H_\bQ)$--trivial part of
$\mathfrak{a}_g(2)$ is one and it is generated by the element
$\om_{12}-\om_{14}$.

\smallskip
{\bf Case (ii)}\qua (the representation $[1^2]$)
\par
\noindent
Let $x_1,\cdots,x_g,y_1,\cdots,y_g$ be the fixed symplectic basis of
$H_\bQ$. Then the element
$x_1\wedge x_2=x_1\otimes x_2-x_2\otimes x_1\in H_\bQ^{\otimes 2}$
is the highest weight vector of the unique summand $[1^2]$
in $H_\bQ^{\otimes 2}$. We define six elements
$\al_{ij}\ (1\leq i<j \leq 4)$ by setting
$$
\al_{ij}=(x_1\wedge x_2)\otimes_{ij} \om_0
$$
where
$$
(u_1\otimes u_2)\otimes_{ij}(u_3\otimes u_4)
$$
is defined to be the element obtained by applying a permutation on
$u_1\otimes u_2\otimes u_3\otimes u_4$
such that $u_1$ and $u_2$ go to the $i$-th and $j$-th places
respectively, while $u_3$ and $u_4$ go to the
$k$-th and $\ell$-th places respectively $(k<\ell)$.
For example
$$
(u_1\otimes u_2)\otimes_{23}(u_3\otimes u_4)
=u_3\otimes u_1\otimes u_2\otimes u_4.
$$
Then it is easy to see that the six elements $\al_{ij}\ (i<j)$
form a basis of $6 [1^2]\subset H_\bQ^{\otimes 4}$.
Now we compute the action of $\sigma_4$ on these elements.
We find that
$$
\al_{12},\ \al_{13},\ \al_{14},\ \al_{23},\ \al_{24},\ \al_{34}
$$
are transformed by $\sigma_4$ to
$$
-\al_{14},\ -\al_{24},\ -\al_{34},\ -\al_{12},\ -\al_{13},\ -\al_{23}.
$$
It follows that the $\bZ/4\bZ$--invariant part is generated by the
following two elements
\begin{equation}
\al^{(1)}=\al_{13}-\al_{24},\quad 
\al^{(2)}=\al_{12}-\al_{14}-\al_{23}+\al_{34}.
\label{eq:11}
\end{equation}

{\bf Case (iii)}\qua (the representation $[2]$)
\par
\noindent
In this case,
the element
$x_1\otimes x_1\in H_\bQ^{\otimes 2}$
is the highest weight vector of the unique summand $[2]$
in $H_\bQ^{\otimes 2}$. We define six elements
$\beta_{ij}\ (1\leq i<j \leq 4)$ by setting
$$
\beta_{ij}=(x_1\otimes x_1)\otimes_{ij} \om_0.
$$
Then it is easy to see that the six elements $\beta_{ij}\ (i<j)$
form a basis of $6 [2]\subset H_\bQ^{\otimes 4}$.
Now we compute the action of $\sigma_4$ on these elements.
We find that
$$
\beta_{12},\ \beta_{13},\ \beta_{14},\ \beta_{23},\ \beta_{24},\ \beta_{34}
$$
are transformed by $\sigma_4$ to
$$
\beta_{14},\ \beta_{24},\ \beta_{34},\ -\beta_{12},\ -\beta_{13},\ -\beta_{23}.
$$
It follows that the $\bZ/4\bZ$--invariant part is generated by the
following single element
$$
\beta_{12}+\beta_{14}-\beta_{23}+\beta_{34}.
$$

{\bf Case (iv)}\qua (the representation $[2^2]$)
\par
\noindent
By \fullref{lem:4}, the multiplicity of
$[2^2]$ in $H_\bQ^{\otimes 4}$ is two.
We have the following two linearly independent elements
\begin{align*}
\ga_1&=(x_1\wedge x_2)\otimes (x_1\wedge x_2)\\
&=(x_1\otimes x_2-x_2\otimes x_1)\otimes (x_1\otimes x_2-x_2\otimes x_1)\\
\ga_2&=x_1\otimes x_1\otimes x_2\otimes x_2+x_2\otimes x_2\otimes 
x_1\otimes x_1\\
&-x_1\otimes x_2\otimes x_1\otimes x_2-x_2\otimes x_1\otimes x_2\otimes x_1
\end{align*}
both of which are highest weight vectors.
The action of $\sigma_4$ is given by
$$
\sigma_4(\ga_1)=-\ga_2, \quad \sigma_4(\ga_2)=-\ga_1.
$$
It follows that $\ga_1-\ga_2$ generates the unique summand
$[2^2]$ which is $\bZ/4\bZ$--invariant.

\smallskip
{\bf Case (v)}\qua (the representation $[31]$)
\par
\noindent
By \fullref{lem:4}, the multiplicity of $[31]$ in $H_\bQ^{\otimes 4}$
is three.
It is easy to see that the element
\begin{align*}
\delta_1=4 x_2\otimes x_1\otimes x_1\otimes x_1-&x_2\otimes x_1\otimes 
x_1\otimes x_1
-x_1\otimes x_2\otimes x_1\otimes x_1\\
-&x_1\otimes x_1\otimes x_2\otimes x_1-x_1\otimes x_1\otimes x_1\otimes x_2
\end{align*}
is the highest weight vector of one particular summand $[31]\subset 
H_\bQ^{\otimes 4}$.
Now define
$$
\delta_2=\sigma_2(\delta_1),\quad
\delta_3=\sigma_3(\delta_1),\quad
\delta_4=\sigma_4(\delta_1).
$$
Here $\sigma_i$ acts on $H_\bQ^{\otimes 4}$ by
$$
\sigma_i(u_1\otimes u_2\otimes u_3\otimes u_4)=
u_i\otimes u_1\otimes\cdots\otimes \hat{u}_i\otimes\cdots\otimes u_4
\quad (i=2,3,4)
$$
(the symbol $\hat{u}_i$ denotes we delete it).
Then the four elements
$$
\delta_1,\delta_2,\delta_3,\delta_4
$$
are transformed by $\sigma_4$ to
$$
\delta_4,\delta_1,\delta_2,\delta_3.
$$
On the other hand, it is easy to see that
$\delta_1+\delta_2+\delta_3+\delta_4=0$ and also that 
$\delta_1,\delta_2,\delta_3$
are linearly independent. Hence there is no summand of
type $[31]$ which is
$\bZ/4\bZ$--invariant.

\smallskip
{\bf Case (vi)}\qua (the representation $[21^2]$)
\par
\noindent
By \fullref{lem:4}, the multiplicity of $[21^2]$ in $H_\bQ^{\otimes 4}$
is three.
It is easy to see that the element
\begin{align*}
\epsilon_1&= x_1\otimes (x_1\wedge x_2\wedge x_3)\\
&=x_1\otimes \sum_{\tau\in\mathfrak{S}_3}
\sgn\tau\, x_{\tau(1)}\otimes x_{\tau(2)}\otimes x_{\tau(3)}
\end{align*}
is the highest weight vector of a certain summand $[21^2]\subset 
H_\bQ^{\otimes 4}$.
Define
$$
\epsilon_2=\sigma_{2} (\epsilon_1),\quad
\epsilon_3=\sigma_{3} (\epsilon_1),\quad
\epsilon_4=\sigma_{4} (\epsilon_1).
$$
Then the four elements
$$
\epsilon_1,\epsilon_2,\epsilon_3,\epsilon_4
$$
are transformed by $\sigma_4$ to
$$
\epsilon_4,\epsilon_1,\epsilon_2,\epsilon_3.
$$
On the other hand, it is easy to check that
$\epsilon_1-\epsilon_2+\epsilon_3-\epsilon_4=0$ and also that
$\epsilon_1,\epsilon_2,\epsilon_3$
are linearly independent.
We can now deduce that the element
$$
\epsilon_1+\epsilon_3=\epsilon_2+\epsilon_4
$$
generates a unique summand of type $[21^2]$ which is
$\bZ/4\bZ$--invariant.

\smallskip
{\bf Case (vii)}\qua (the representation $[1^4]$)
\par
\noindent
In this case, there is a unique $[1^4]$ in
$H_\bQ^{\otimes 4}$ whose highest weight vector is
$$
x_1\wedge x_2\wedge x_3\wedge x_4
=\sum_{\tau\in\mathfrak{S}_4} \sgn \tau\, x_{\tau(1)}\otimes
x_{\tau(2)}\otimes x_{\tau(3)}\otimes
x_{\tau(4)}.
$$
This element is not cyclically invariant. Hence $[1^4]$ does not appear in
$\mathfrak{a}_g(2)$. In general, $[1^{k+2}]$ appears in
$\mathfrak{a}_g(k)$ if and only if $k$ is odd.

\smallskip
{\bf Case (viii)}\qua (the representation $[4]$)
\par
\noindent
In this case, clearly $[4]=S^4 H_\bQ$ is $\bZ/4\bZ$--invariant so that
it appears in $\mathfrak{a}_g(2)$. In general, $[k+2]=S^{k+2} H_\bQ$
appears in $\mathfrak{a}_g(k)$ for any $k$.

This completes the proof.
\end{proof}

Now we prove the main theorem of this paper.

\proof[Proof of \fullref{th:main}]
Our task is to prove the exactness of the
sequence
\begin{equation}
\Lambda^2 \mathfrak{a}_{g}(1)\
\overset{[\ , \ ]}{\lra}\
\mathfrak{a}_{g}(2)\
\overset{C_{13}}{\lra}\
\Lambda^2 H_\bQ/\bQ(\om_0)\cong [1^2]\
\lra\ 0
\label{eq:sq}
\end{equation}
which implies the required isomorphism
$$
H_1(\mathfrak{a}^+_{g})_2 \cong \Lambda^2 H_\bQ/\bQ(\om_0).
$$
We begin by showing that the homomorphism
$C_{13}$ vanishes identically on commutators.
For this, we consider $\mathfrak{a}_{g}(1)$ to be a submodule
of $H_\bQ^{\otimes 3}$ by \fullref{prop:al}.
Let
\begin{align*}
\xi&=\sum_i u^i_1\otimes u^i_2\otimes u^i_3\in \mathfrak{a}_{g}(1)
\subset H_\bQ^{\otimes 3}\\
\eta&=\sum_j v^j_1\otimes v^j_2\otimes v^j_3\in \mathfrak{a}_{g}(1)
\subset H_\bQ^{\otimes 3}
\end{align*}
be any two elements of $\mathfrak{a}_{g}(1)$.
Then the bracket $[\xi,\eta]\in \mathfrak{a}_{g}(2)$
is given by $[\xi,\eta]=\sum_{i,j}[\xi,\eta]_{ij}$ where
\begin{align*}
[\xi&,\eta]_{ij}\\
=&(u^i_1\cdot v^j_2)\, v^j_1\otimes u^i_2\otimes u^i_3\otimes v^j_3
+(u^i_1\cdot v^j_3)\, v^j_1\otimes v^j_2\otimes u^i_2\otimes u^i_3\\
-&(v^j_1\cdot u^i_2)\, u^i_1\otimes v^j_2\otimes v^j_3\otimes u^i_3
-(v^j_1\cdot u^i_3)\, u^i_1\otimes u^i_2\otimes v^j_2\otimes v^j_3.
\end{align*}
It follows that
\begin{align*}
C_{13}&([\xi,\eta]_{ij})\\
=&(u^i_1\cdot v^j_2)(v^j_1\cdot u^i_3) \, u^i_2\otimes v^j_3
+(u^i_1\cdot v^j_3)(v^j_1\cdot u^i_2)\, v^j_2\otimes u^i_3\\
-&(v^j_1\cdot u^i_2)(u^i_1\cdot v^j_3) \, v^j_2\otimes u^i_3
-(v^j_1\cdot u^i_3)(u^i_1\cdot v^j_2) \, u^i_2\otimes v^j_3\\
=& 0
\end{align*}
and hence $C_{13}([\xi,\eta])=0$ as required.
We will later observe that the vanishing of $C_{13}$ on commutators
is valid in a broader context (see \fullref{prop:c}).

Next we compute the value of $C_{13}$ on
the component of type $[1^2]$ in $\mathfrak{a}_{g}(2)$.
In the notation of Case (ii) of \fullref{prop:a2}, we have
$$
\al_{12}=\sum_{i=1}^g (x_1\otimes x_2 -x_2\otimes x_1)\otimes
  (x_i\otimes y_i-y_i\otimes x_i)
$$
so that
$$
C_{13}(\al_{12})=-x_2\otimes x_1+x_1\otimes x_2.
$$
Similarly, we have
\begin{align*}
&C_{13}(\al_{13})=0,\\
&C_{13}(\al_{14})=x_1\otimes x_2-x_2\otimes x_1,\quad
C_{13}(\al_{23})=x_1\otimes x_2-x_2\otimes x_1,\\
&C_{13}(\al_{24})=2g\, (x_1\otimes x_2-x_2\otimes x_1),\quad
C_{13}(\al_{34})=x_1\otimes x_2-x_2\otimes x_1.
\end{align*}
It follows that the values of $C_{13}$ on
the two elements $\al^{(1)}, \al^{(2)}$
(see \eqref{eq:11} in \fullref{prop:a2})
are given by
$$
C_{13}(\al^{(1)})=-2g\, (x_1\otimes x_2-x_2\otimes x_1),\quad
C_{13}(\al^{(2)})=0.
$$
Thus we have proved the exactness at the last factor of the
sequence \eqref{eq:sq},
namely the surjectivity of the homomorphism
$$
C_{13}\co\mathfrak{a}_{g}(2)\lra \Lambda^2 H_\bQ/\bQ(\om_0).
$$
It remains to prove that all the irreducible components appearing
in \fullref{prop:a2} other than the one generated by
the element $\al^{(1)}$ above can be represented as commutators.
Consider the degree two part
${\mathfrak l}_g(2)$ of the
graded Lie algebra $\mathfrak{h}^\bQ_{g,1}$
which is a graded Lie subalgebra of $\mathfrak{a}^+_{g}$.
It is a consequence of a fundamental paper \cite{Hain97} by Hain,
suitably adapted to the case with one boundary component as in \cite{Morita96},
that the bracket operation
$$
\Lambda^2 {\mathfrak l}_g(1)\
\overset{[\ .\ ]}{\lra}\
{\mathfrak l}_g(2)\cong \bQ\oplus [1^2]\oplus [2^2]
$$
is surjective (see also \cite{Morita01} for a related result).

Thus we have only to prove that the three components
$[2], [21^2], [4]$ appearing in the decomposition of
$\mathfrak{a}_{g}(2)$ belong to the commutator ideal
of $\mathfrak{a}^+_{g}$. Recall from \fullref{ex:1}
that $\mathfrak{a}_{g}(1)\cong \Lambda^3 H_\bQ\oplus S^3 H_\bQ$.
We consider case by case.

\smallskip
{\bf Case of the representation $[2]$}\qua
We consider two elements
\begin{align*}
\xi_1&= x_2\wedge y_2\wedge x_3\in \Lambda^3 H_\bQ\subset \mathfrak{a}_{g}(1)\\
\eta_1&= x_1\otimes x_1\otimes y_3+x_1\otimes y_3\otimes x_1
+y_3\otimes x_1\otimes x_1
\in S^3 H_\bQ\subset \mathfrak{a}_{g}(1).
\end{align*}
Then explicit computation implies
\begin{align*}
&[\xi_1, \eta_1]\\
=&-x_1\otimes x_1\otimes (x_2\otimes y_2-y_2\otimes x_2)
-(x_2\otimes y_2-y_2\otimes x_2)\otimes x_1\otimes x_1\\
&-x_1\otimes (x_2\otimes y_2-y_2\otimes x_2)\otimes x_1\\
&+x_2\otimes x_1\otimes x_1\otimes y_2-y_2\otimes x_1\otimes x_1\otimes x_2.
\end{align*}
If we apply the contraction
$C_{11}\co H_\bQ^{\otimes 4}\ra H_\bQ^{\otimes 2}$ defined by
$$
C_{11}(u_1\otimes u_2\otimes u_3\otimes u_4)=(u_1\cdot u_2)\, u_3\otimes u_4
$$
to the above element $[\xi_1, \eta_1]$, we obtain
$$
-2 x_1\otimes x_1
$$
which is the highest weight vector of $[2]\subset H_\bQ^{\otimes 2}$.
It follows that the summand $[2]$ is a commutator.

\smallskip
{\bf Case of the representation $[21^2]$}\qua
We consider two elements
\begin{align*}
\xi_2&= x_1\wedge x_2\wedge x_3\in \Lambda^3 H_\bQ\subset \mathfrak{a}_{g}(1)\\
\eta_2&= x_1\otimes x_1\otimes y_1+x_1\otimes y_1\otimes x_1
+y_1\otimes x_1\otimes x_1
\in S^3 H_\bQ\subset \mathfrak{a}_{g}(1).
\end{align*}
Then explicit computation implies
\begin{align*}
&[\xi_2, \eta_2]\\
=&x_1\otimes \left(x_1\otimes (x_2\otimes x_3-x_3\otimes x_2)-
(x_2\otimes x_3-x_3\otimes x_2)\otimes x_1\right)\\
-& x_2\otimes (x_1\otimes x_1\otimes x_3-x_3\otimes x_1\otimes x_1)\\
+& x_3\otimes (x_1\otimes x_1\otimes x_2-x_2\otimes x_1\otimes x_1).
\end{align*}
This is the highest weight vector of the summand $[21^2]\subset 
\mathfrak{a}_{g}(2)$
because if we replace $x_3$ with $x_2$ in the above tensor, we obtain
$0$ and the same is true if we replace $x_2$ with $x_1$.

\smallskip
{\bf Case of the representation $[4]$}\qua
We consider the element
$$
\xi_3= x_1\otimes x_1\otimes x_1\in S^3 H_\bQ\subset \mathfrak{a}_{g}(1)
$$
and the above element $\eta_2\in S^3 H_\bQ\subset \mathfrak{a}_{g}(1)$.
Then explicit computation shows
$$
[\xi_3, \eta_2]
=4 x_1\otimes x_1\otimes x_1\otimes x_1
$$
which is the highest weight vector of the summand $[4]\subset 
\mathfrak{a}_{g}(2)$.

This completes the proof of \fullref{th:main}.\endproof

\section{Cycles on the moduli spaces}

There are several method of constructing explicit
(co)cycles of the mapping class group as well as the moduli space of curves.
They include group cohomological construction given in
\cite{Morita96}, combinatorial construction using the natural cell
structures on the moduli spaces by Witten and Kontsevich, Kontsevich's
construction \cite{Kontsevich93,Kontsevich94}
using $A_\infty$ algebras with scalar products,
and more recently Penner and the author \cite{MP}
obtained another combinatorial method.
However all the classes obtained in these ways turn out to be
expressed as polynomials in the Mumford--Morita--Miller tautological classes
as shown in the work of Kawazumi and the author \cite{KM96,KM01}
for the group cohomological cocycles
and Igusa \cite{Igusa} and Mondello \cite{Mondello}
for the Witten and Kontsevich cycles.

Recently Madsen and Weiss \cite{MW} proved that the
stable cohomology of the mapping class group is the polynomial
algebra generated by the Mumford--Morita--Miller classes.
Hence we can say that there are many ways of constructing cocycles
for any stable cohomology class. On the other hand, the works of
Harer and Zagier \cite{HZ} and Penner \cite{Penner88}
determining the orbifold Euler characteristics of the moduli space
imply that there should exist huge amount of unstable classes.
However there are very few works on the unstable classes
(see Looijenga \cite{Looijenga93} and Tommasi \cite{Tommasi}).
It would be worthwhile to try to use the following
(suitably modified) theorem of Kontsevich
to produce unstable homology classes of the moduli
space.

\begin{theorem}[Kontsevich \cite{Kontsevich93,Kontsevich94}]\label{th:Kontsevich}
There is an isomorphism
$$
PH^k(\mathfrak{a}^+_\infty)^{\mathrm{Sp}}_{2n}\cong
\bigoplus_{2g-2+m=n>0,\, m>0} H_{2n-k}(\mathbf{M}_g^m;\bQ)_{\mathfrak{S}_m}.
$$
\end{theorem}

Here $P$ denotes the {\it primitive part} of
$H^*(\mathfrak{a}^+_\infty)^{\mathrm{Sp}}$ which has a
natural structure of Hopf algebra, the subscript
$2n$ denotes the weight $2n$ part and
$\mathbf{M}_g^m$ denotes the moduli space of genus
$g$ smooth curves with $m$ punctures.

If we apply the above theorem to Tommasi's unstable cohomology class
in $H^5(\mathbf{M}_4;\bQ)$, for example, it should yield certain
elements in $PH^{k}(\mathfrak{a}^+_\infty)^{\mathrm{Sp}}_{14}$
for $k=7,9$.

Now if we consider the abelianization
$
\mathfrak{a}_{g}^+\ra H_1(\mathfrak{a}_{g}^+)
$
of the Lie algebra $\mathfrak{a}_{g}^+$ and apply the above
theorem of Kontsevich, we obtain a homomorphism
$$
\lim_{g\to\infty}
PH^k\left(H_1(\mathfrak{a}^+_g)\right)^{\mathrm{Sp}}_{2n}\lra
\bigoplus_{2g-2+m=n>0,\, m>0} H_{2n-k}(\mathbf{M}_g^m;\bQ)_{\mathfrak{S}_m}.
$$
We can then ask which part of the cohomology of the moduli spaces
can be obtained in this way.
In particular, if the abelianization of $\mathfrak{a}^+_g$
was large like in the case of ${\mathfrak l}_g^+$, then we would obtain
many candidates for unstable classes of the moduli spaces.

Now our main result (\fullref{th:main}) together with \fullref{ex:1}
gives rise to a surjective homomorphism
$$
\mathfrak{a}^+_g\lra S^3H_\bQ\oplus\Lambda^3 H_\bQ
\oplus\Lambda^2 H_\bQ/\bQ(\om_0)
$$
and hence
\begin{align*}
\lim_{g\to\infty} PH^{k}&(S^3 H_\bQ\oplus\Lambda^3 H_\bQ\oplus
\Lambda^2 H_\bQ/\bQ(\om_0))^{\mathrm{Sp}}\\
&\lra \bigoplus_{2g-2+m=n>0,\, m>0} 
H_{2n-k}(\mathbf{M}_g^m;\bQ)_{\mathfrak{S}_m}.
\end{align*}
It can be shown that the source of the above homomorphism contain many
cohomology classes.
By restricting further to the weight $2$ part, we obtain
$$
\lim_{g\to\infty} PH^{k}(\Lambda^2 H_\bQ/\bQ(\om_0))^{\mathrm{Sp}}
\lra
\bigoplus_{k=2g-2+m>0,\, m>0} H_{k}(\mathbf{M}_g^m;\bQ)_{\mathfrak{S}_m}
$$
whose left hand side can be determined as follows.

\begin{proposition}\label{prop:sp}
The $\mathrm{Sp}$--trivial part
$H^*\left(\Lambda^2 H_\bQ/\bQ(\om_0)\right)^{\mathrm{Sp}}$ of the 
cohomology group
of the irreducible representation $\Lambda^2 H_\bQ/\bQ(\om_0)\cong [1^2]$
is isomorphic to the cohomology of the product
$
S^{5}\times S^{9}\times S^{13}\times \cdots
$
as $g$ tends to $\infty$ so that we have
$$
\lim_{g\to\infty} PH^k\left(\Lambda^2 H_\bQ/\bQ(\om_0)\right)^{\mathrm{Sp}}
\cong
\begin{cases}
\bQ &\quad (k\equiv 1\ \mathrm{mod}\ 4, k\not= 1)\\
0 &\quad (\mathrm{otherwise}).
\end{cases}
$$
\end{proposition}

\begin{proof}
By the classical invariant theory of Weyl, any $\mathrm{Sp}$--invariant
cohomology class in degree $k$ of the $\mathrm{Sp}$--module
$\Lambda^2 H_\bQ/\bQ(\om_0)$
can be expressed as a
linear combination of certain contractions
$$
C\co  \Lambda^k(\Lambda^2 H_\bQ/\bQ(\om_0))\ra \bQ
$$
each of which is a
various multiples of the intersection pairing on $H_\bQ$.
Such contractions can be enumerated by bivalent graphs and primitivity
in cohomology corresponds to connectedness of graphs.
There is only one connected bivalent graph with $k$ vertices,
namely a $k$--gon, for each $k\geq 2$.
Let $C_k$ be the corresponding contraction.
Now we apply $C_k$
on the following particular element
$$
\lambda_k=\sum_{\tau\in\mathfrak{S}_k}
\sgn\tau\, \lambda_k^{(\tau(1))}\otimes\cdots\otimes\lambda_k^{(\tau(k))}
$$
of $\Lambda^k(\Lambda^2 H_\bQ/\bQ(\om_0))$,
where
\begin{align*}
&\lambda_k^{(1)}=x_1\wedge y_2, \
\lambda_k^{(2)}=x_2\wedge y_3,\
\cdots, \\
&\lambda_k^{(k-1)}=x_{k-1}\wedge y_k,\
\lambda_k^{(k)}=x_k\wedge y_1.
\end{align*}
By explicit computations, we see that
$$
C_k(\lambda_k)\not=0
$$
if and only if $k\equiv 1\mod 4\ (k\not=1)$.
On the other hand, it can be checked easily that
any contraction corresponding to {\it disconnected}
graph vanishes on $\lambda_k$.
The result follows from these facts.
\end{proof}

\begin{remark}
If we replace $\Lambda^2 H_\bQ/\bQ(\om_0)$ by the second symmetric power
$S^2H_\bQ$ in the above argument, we obtain similar result
where $k\equiv 1\ \mathrm{mod}\ 4$ is replaced by
$k\equiv 3\ \mathrm{mod}\ 4$. This can be also obtained by making
use of the fact that $S^2H_\bQ$ is isomorphic to the Lie algebra
of $\mathrm{Sp}(2g,\bQ)$.
\end{remark}

Thus we obtain infinitely many unstable
{\it odd} dimensional homology classes of the
moduli spaces.

\begin{remark}
If we apply \fullref{th:Kontsevich}
to \fullref{prop:sp}, we obtain homomorphisms
$$
\lim_{g\to\infty} PH^k\left(\Lambda^2 H_\bQ/\bQ(\om_0)\right)^{\mathrm{Sp}}
\lra H_{4k+1}(\mathbf{M}_1^{4k+1})_{\mathfrak{S}_{4k+1}}
$$
for all $k\geq 1$. Recently, Conant \cite{Conant} proved
that the cycles of the genus one moduli spaces obtained in this way are all
non-trivial. It would be interesting to compare this result with
that of Getzler \cite{Getzler} where the rational cohomology
of genus one moduli spaces are determined.

\end{remark}

\section{Further results}

It seems worthwhile to consider the Lie algebra consisting of derivations
on free associative algebras {\it without} the symplectic constraint.
For this, let $H_n$ denote a vector space over $\bQ$ of
dimension $n\geq 2$.
Also let $T(H_n)$ denote the tensor algebra without unit
generated by $H_n$. Now define
$$
\mathrm{Der}(T(H_n))=\{\mathrm{derivation}\ D\ \mathrm{of}\
\mathrm{the}\ \mathrm{free}\ \mathrm{associative}\ \mathrm{algebra}\
T(H_n)\}.
$$
It has a natural structure of a graded Lie algebra
$$
\mathrm{Der}(T(H_n))=\bigoplus_{k=0}^\infty
\mathrm{Der}(T(H_n))(k)
$$
$$
\mathrm{Der}(T(H_n)(k)=
\Hom \left(H_n, H_n^{\otimes (k+1)}\right)
\cong H^*_n\otimes H_n^{\otimes (k+1)}.\leqno{\rm where}
$$
$$
\mathrm{Der}^+(T(H_n))=\bigoplus_{k=1}^\infty
H^*_n\otimes H_n^{\otimes (k+1)}\leqno{\rm Let}
$$
be the ideal of
$\mathrm{Der}(T(H_n))$ consisting of derivations with
positive degrees. By definition,
$\mathfrak{a}^+_{g}$ is clearly a Lie subalgebra of
$\mathrm{Der}^+(T(H_{2g}))$ and we are interested in the
structure of $\mathrm{Der}^+(T(H_n))$ as well.
In particular, we have the following problem.

\begin{problem}
Determine whether the Lie algebra $\mathrm{Der}^+(T(H_n))$ is finitely 
generated
or not. In particular, is the abelianization
$$
H_1(\mathrm{Der}^+(T(H_n)))
$$
of $\mathrm{Der}^+(T(H_n))$
finitely generated or not?
\end{problem}

\begin{remark}
We can also consider the free Lie algebra $\mathcal{L}(H_n)$ generated by
$H_n$ and also the Lie algebra $\mathrm{Der}(\mathcal{L}(H_n))$ consisting
of derivations of it. This Lie algebra can be naturally considered as
a Lie subalgebra of $\mathrm{Der}(T(H_n))$. In this case, however,
we know that $\mathrm{Der}^+(\mathcal{L}(H_n))$ is {\it not} finitely
generated because the surjective homomorphisms
$\mathrm{trace}(k)\co\mathrm{Der}(\mathcal{L}(H_n))(k)\ra S^k H_n$
defined in \cite{Morita93a}
vanish on the commutator ideal
$[\mathrm{Der}^+(\mathcal{L}(H_n)),\mathrm{Der}^+(\mathcal{L}(H_n))]$
for any $k\geq 2$
(see \cite{Morita06}).
\end{remark}

The Lie algebra $\mathrm{Der}^+(T(H_n))$ is graded. It follows that
the abelianization\break $H_1(\mathrm{Der}^+(T(H_n)))$ is bigraded
and we have a direct sum decomposition
$$
H_1(\mathrm{Der}^+(T(H_n)))=
\bigoplus_{m=1}^\infty
H_1(\mathrm{Der}^+(T(H_n)))_m.
$$
Clearly we have
$$
H_1(\mathrm{Der}^+(T(H_n)))_1=H^*_n\otimes H_n^{\otimes 2}.
$$
Now we show that the next summand, namely the weight $2$ part of the
abelianization of $\mathrm{Der}^+(T(H_n))$ is non-trivial.
For this, consider the contraction
$$
C_{13}\co  \mathrm{Der}(T(H_n))(2)\cong
H^*_n\otimes H_n^{\otimes 3}\lra H_n^{\otimes 2}
$$
defined by
$$
C_{13}(f\otimes u_1\otimes u_2\otimes u_3)=f(u_2)\, u_1\otimes u_3
$$
where $f\in H^*_n, u_i\in H_n$. It is clearly surjective.

\begin{proposition}\label{prop:c}
The composition
$$
\Lambda^2\mathrm{Der}(T(H_n))(1)\
\overset{[\ ,\ ]}\lra\
\mathrm{Der}(T(H_n))(2)\
\overset{C_{13}}\lra\
H_n^{\otimes 2}
$$
is trivial. In other words, the homomorphism $C_{13}$ is trivial
on the commutator ideal
$[\mathrm{Der}^+(T(H_n)),\mathrm{Der}^+(T(H_n))]$
of the Lie algebra $\mathrm{Der}^+(T(H_n))$ so that
it induces a surjection
$$
H_1(\mathrm{Der}^+(T(H_n)))_2\lra H_n^{\otimes 2}.
$$
\end{proposition}

\begin{proof}
The bracket operation
on $\Lambda^2\mathrm{Der}(T(H_n))(1)$ is given by
\begin{align*}
[f\otimes u_1\otimes& u_2, h\otimes v_1\otimes v_2]\\
=\ & f(v_1)\, h\otimes u_1\otimes u_2\otimes v_2
+ f(v_2)\, h\otimes v_1\otimes u_1\otimes u_2\\
  -& h(u_1)\, f\otimes v_1\otimes v_2\otimes u_2
-h(u_2)\, f\otimes u_1\otimes v_1\otimes v_2
\end{align*}
where $f,h\in H^*_n, u_i, v_i\in H_n$.
It follows that
\begin{align*}
C_{13}([f\otimes u_1&\otimes u_2, h\otimes v_1\otimes v_2])\\
=\ & f(v_1)h(u_2)\, u_1\otimes v_2+ f(v_2)h(u_1)\, v_1\otimes u_2\\
  -& h(u_1)f(v_2)\, v_1\otimes u_2-h(u_2)f(v_1)\, u_1\otimes v_2\\
=\ &0.
\end{align*}
This completes the proof.
\end{proof}

\begin{remark}
There is a close connection between the contraction $C_{13}$ and the
trace map. In fact, the equality
$C_{13}=-2\, \mathrm{trace}(2)$ holds so that
the following diagram is
commutative (up to a non-zero scalar)
\begin{equation*}
\begin{CD}
\mathrm{Der}(\mathcal{L}(H_n))(2) @>{\mathrm{trace(2)}}>> S^2 H_n \\
@VVV @VVV \\
\mathrm{Der}(T(H_n))(2) @>{C_{13}}>> H_n^{\otimes 2}.
\end{CD}
\end{equation*}
\end{remark}

\begin{conjecture}
The contraction $C_{13}\co \mathrm{Der}(T(H_n))(2)\ra
H_n^{\otimes 2}$ induces an isomorphism
$$
H_1(\mathrm{Der}^+(T(H_n)))_2\cong H_n^{\otimes 2}.
$$
\end{conjecture}

By an explicit computation using a computer, we have obtained
the following result.

\begin{proposition}
The above conjecture is true for the case $n=2$, namely
we have an isomorphism
$$
H_1(\mathrm{Der}^+(T(H_2)))_2\cong H_2^{\otimes 2}.
$$
Furthermore, the homomorphism
$$
\mathrm{Der}(T(H_2))(1)\otimes \mathrm{Der}(T(H_2))(k)
\ \overset{[\ , \ ]}{\lra}\
  \mathrm{Der}(T(H_2))(k+1)
$$
induced by the bracket operation
is surjective for $k=2,3$ so that
$$
H_1(\mathrm{Der}^+(T(H_2)))_m=0
$$
for $m=3,4$.
\end{proposition}

\begin{remark}
At present, we have no idea about the (non)-triviality of the
abelianization $H_1(\mathrm{Der}^+(T(H_n)))_m$ for $n\geq 3$ and
$m\geq 3$.  In the case of $\mathrm{Der}^+(\mathcal{L}(H_n))$,
Kassabov \cite{Kassabov} has determined the abelianization of it in a
certain stable range (cf \cite{Morita06}).
\end{remark}

\bibliographystyle{gtart}
\bibliography{link}

\end{document}